\documentclass{amsart}

\usepackage[english]{babel}
\usepackage[utf8]{inputenc}
\usepackage{amsmath}
\usepackage{amsthm}
\usepackage{amsfonts}
\usepackage{amssymb}
\usepackage{enumitem}

\theoremstyle{definition}
\newtheorem{defin}{Definition}

\newtheorem{teo}{Theorem}

\newtheorem{prop}{Proposition}

\theoremstyle{remark}
\newtheorem{rem}{Remark}[section]

\begin{document}
\title{Newton non-degenerate foliations and blowing-ups}
\author{Beatriz Molina-Samper}
\address{Dpto. \'Algebra, An\'alisis Matem\'atico, Geometr\'ia y Topolog\'ia and Instituto de Matem\'aticas de la Universidad de Valladolid, Facultad de Ciencias, Universidad de Valladolid, Paseo de Bel\'en, 7, 47011 Valladolid, Spain}
\email{beatriz.molina@uva.es}
\thanks{Supported by the MINECO of Spain, with the MTM2016-77642-C2-1-P research project and by the MECD of Spain with the FPU14/02653 grant.}			
	
\date{}
\dedicatory{}

\begin{abstract}
A codimension one singular holomorphic foliation is Newton non-degenerate if it satisfies some non-degeneracy conditions, in terms of its Newton polyhedra system. These conditions are similar to the ones of Kouchnirenko and Oka for the case of functions. We introduce the concept of logarithmic reduction of singularities and we prove that a foliation is Newton non-degenerate if and only if it admits a logarithmic reduction of singularities of a combinatorial nature.
\end{abstract}

\maketitle
	
\section{Introduction}
A foliation $\mathcal{F}$ on a complex space $M$ admits a combinatorial logarithmic reduction of singularities with respect to a normal crossings divisor $E \subset M$ when the centers of blowing-up are compatible with the natural stratification provided by $E$. The data $(M,E; \mathcal{F})$ is called a foliated space. The goal of this paper is to characterize the class of codimension one foliations admitting such reductions of singularities. We give this characterization in terms of the Newton polyhedra associated to the pair $(\mathcal{F},E)$, and we state the Equivalence Theorem as follows:

\begin{quote}
\textbf{Theorem \ref{teo:equivalencia}} A foliated space admits a combinatorial logarithmic reduction of singularities if and only if it is Newton non-degenerate.
\end{quote}

A foliated space $(M,E;\mathcal{F})$ is logarithmically desingularized, or logarithmically regular, at a point $p\in M$ if there is a local generator $\eta$ of $\mathcal{F}$ at $p$ that can be written as
$$
\eta=\sum_{j=1}^ea_j \frac{dx_j}{x_j}+\sum_{j=e+1}^n a_jdy_j,
$$
where $E=(\prod_{j=1}^ex_j=0)$ and there is at least one unit among the coefficients $a_1,a_2,\ldots,a_n$. A combinatorial logarithmic reduction of singularities intends to transform the foliated space into a logarithmically regular one.

This kind of reduction of singularities is closely related with the classical one for foliations, whose development starts at the work of Seidenberg \cite{Sei}, in dimension two, and at the papers \cite{Can,Can-C}, in higher dimension. More precisely, when there are no ``hidden saddle-nodes'' in the foliation, the logarithmically regular points coincide with the classical presimple points in the cited works.

The concept of Newton non-degenerate foliated space is given by extending the classical ideas for germs of hypersurface in the works \cite{Kou,Oka}. We attach a polyhedron $N_S$ to each stratum $S$ of the natural stratification induced by $E$. Each compact face of $N_S$ provides a weighted initial form for a local  logarithmic generator of the foliation. We ask the weighted initial forms to have no zeros in the corresponding spaces of the form $(S,p)\times T$, where $T$ is a complex torus, at any point $p\in S$. When these non-degeneracy conditions hold, we say that the foliated space is Newton non-degenerate.

The proof of the Equivalence Theorem is based in two fundamental results. The first one is the stability of being Newton non-degenerate under combinatorial blowing-ups and blowing-downs. The second one is the equivalence between being Newton non-degenerate and being logarithmically desingularized, under the hypothesis of having a desingularized polyhedra system. In this way, we conclude by applying a reduction of singularities for the polyhedra system, whose existence is proved in \cite{Mol}.

Frobenius integrability condition does not intervene in the theorem. Hence, we have a similar statement for fields of hyperplanes, although the geometrical interpretation is different than in the case of codimension one foliations.

\section{Combinatorial Logarithmic Reduction of Singularities}
Let $M$ be a $n$-dimensional nonsingular complex analytic space. A \emph{foliation} $\mathcal{F}$ on $M$ is an invertible coherent $\mathcal{O}_M$-submodule $\mathcal{F}\subset \Omega_M^1$, integrable and saturated in the sense that $\mathcal{F}=\mathcal{F}^{\bot\bot}$. In local terms, a foliation is generated by a Frobenius integrable germ of holomorphic one-form
$$
\omega=f_1dz_1+f_2dz_2+\cdots+f_ndz_n,
$$
without common factors in its coefficients.

For our purposes, we need to consider normal crossings divisors on $M$. A \emph{normal crossings divisor} $E=\cup_{i\in I} E_i$ is a finite union of smooth hypersurfaces $E_i$ such that $E$ can be seen as a union of coordinate hyperplanes, locally at each point $p\in M$. Besides to this usual definition, we ask to the sets $E_J=\cap_{j\in J}E_j$ to be connected for any $J\subset I$; that is the $E_J$ are the adherence of the strata defined by $E$. We refer to this additional condition by saying that $E$ is a \emph{strong normal crossings divisor}. We say that such a pair $(M,E)$ is an \emph{ambient space}. We define a \emph{combinatorial blowing-up between ambient spaces} as a map
$$
\pi:(M',E')\to (M,E), \quad E'=\pi^{-1}(E)
$$
induced by a blowing-up $M' \to M$ centered at one of the sets $E_J\subset E$.

A \emph{foliated space $(M,E;\mathcal{F})$} is the datum of an ambient space $(M,E)$ and a foliation $\mathcal{F}$ on $M$. A combinatorial blowing-up $\pi:(M',E')\to (M,E)$ is \emph{admissible for $(M,E;\mathcal{F})$} when the center $E_J$ is invariant for $\mathcal{F}$. We write, for short
$$
\pi:(M',E';\mathcal{F}')\to (M,E; \mathcal{F}),
$$
where $\mathcal{F}'$ is the transform of $\mathcal{F}$ by $\pi$.

We consider the set $\mathcal{H}_{M,E}$ of the subsets $J \subset I$, such that $E_J\ne \emptyset$. Given $J\in \mathcal{H}_{M,E}$, the stratum $S_J$ is 
$$
S_J=E_J\setminus \cup_{j \in I\setminus J}E_j.
$$
In order to give labels in a convenient manner for local coordinate systems at the points of $S_J$, we make a choice of a set $c(J)$, for each $J\in \mathcal{H}_{M,E}$, such that $J\cap c(J)=\emptyset$ and $\#J \cup c(J)=n$. For each point $p\in S_J$, a \emph{local coordinate system $(\mathbf{x},\mathbf{y})$ adapted to $E$} is a pair of families of germs of functions
$$
\mathbf{x}=(x_j)_{j\in J},\;\mathbf{y}=(y_j)_{j\in c(J)},
$$
such that their union forms a local coordinate system at $p$ and the divisor $E$ is given by $E=(\prod_{j\in J}x_j=0)$. 

Taking a logarithmic point of view with respect to the divisor $E$, we consider $\mathcal{F}$ as being locally defined at a point $p\in M$ by a meromorphic one-form
\begin{equation} \label{eq:logaritmico}
\eta =\sum_{j\in J}a_j\frac{dx_j}{x_j}+\sum_{j\in c(J)}a_jdy_j,
\end{equation}
without common factors in its coefficients, where $J$ corresponds to the stratum $S_J$ such that $p\in S_J$.  This expression allows us to define the \emph{logarithmic singular locus $\operatorname{logSing}(\mathcal{F},E)$} as follows
$$
\text{logSing}(\mathcal{F},E)=\{p \in M ; \; \nu_p(\mathcal{F},E) > 0\},
$$
where $\nu_p(\mathcal{F},E)$ is the minimum of the orders at $p$ of the coefficients of $\eta$.

We say that an $E_J$ is a \emph{log-admissible center of blowing-up for $(M,E;\mathcal{F})$} when $E_J\subset \text{logSing}(\mathcal{F},E)$.

\begin{rem}
Note that if $E_J$ is log-admissible, it is also admissible.
\end{rem}
\begin{defin}
A foliated space $(M,E;\mathcal{F})$ is of \emph{logarithmic toric type} if there is a finite sequence of log-admissible combinatorial blowing-ups
\begin{equation}  \label{eq:logadmisible}
(M',E';\mathcal{F}') \rightarrow \cdots \rightarrow (M_1,E^1;\mathcal{F}_1) \rightarrow (M,E;\mathcal{F})
\end{equation}
such that $\text{logSing}(\mathcal{F}',E')=\emptyset$.
\end{defin}

\begin{rem} \label{rem:casoCH}
The reduction of singularities for holomorphic foliations \cite{Can,Can-C,Sei}, attempts two objectives: to obtain either presimple points or, the more restrictive, simple points. Roughly speaking, simple points are ``presimple ones without resonances''.  There is a context in which presimple points coincide with logarithmically non-singular points: the case of complex hyperbolic foliations (see \cite{Mol2}).  We recall that a foliation $\mathcal{F}$ on $M$ is complex hyperbolic (see \cite{Can-R}) if there is no holomorphic map $\phi:(\mathbb{C}^2,0) \rightarrow M$ such that $0$ is a saddle-node for $\phi^{-1}\mathcal{F}$. In the two dimensional case, being complex hyperbolic is equivalent to have a reduction of singularities without saddle-nodes; this is the case considered in \cite{Cam-N-S}.

Hence, in the complex hyperbolic context, a foliated space is of logarithmic toric type if and only if it has a ``combinatorial reduction of singularities to presimple points'' (weak toric type in \cite{Mol3}). Then, as a direct consequence of Theorem \ref{teo:equivalencia}, we have that a complex hyperbolic foliated space is Newton non-degenerate if and only if it is of weak toric type.
\end{rem}

\section{Newton Non-degenerate Foliated Spaces}
We devote this section to introduce the definition of Newton non-degenerate foliated spaces $(M,E;\mathcal{F})$. In order to do it, we use Newton polyhedra systems, following the definitions introduced in \cite{Mol}.

The set $\mathcal{H}_{M,E}$, introduced before, is called \emph{support fabric} in \cite{Mol}; as we have seen, it represents the natural stratification on $M$ induced by $E$. We associate to each $J\in \mathcal{H}_{M,E}$ a positively convex polyhedron $N_J \subset \mathbb{R}_{\geq 0}^J$ as follows. Take a point $p$ in the stratum $S_J$ and a local logarithmic generator $\eta$ of $\mathcal{F}$ at $p$ as in Equation \ref{eq:logaritmico}. We decompose each coefficient $a_j$ of $\eta$ as
\begin{equation}\label{eq:decompsigma}
a_j=\sum_{\sigma\in \mathbb{Z}_{\geq 0}^J} a_{j,\sigma} (\mathbf{y}) \mathbf{x}^\sigma, \quad \mathbf{x}^\sigma=\mbox{$\prod\limits_{j\in J}$}x_j^{\sigma(j)}.
\end{equation}
The polyhedron $N_J$ is the positively convex hull in $\mathbb{R}^J$ of the set $\text{Supp}(\eta; (\mathbf{x},\mathbf{y}))$, where
$$
\text{Supp}(\eta; (\mathbf{x},\mathbf{y}))=\{\sigma \in \mathbb{Z}^J; \; \text{ there is } j\in J\cup c(J) \text{ with } a_{j,\sigma}(\mathbf{y}) \neq 0\}.
$$
\begin{rem}
The definition of $N_J$ is independent of the choice of $p\in S_J$, the particular adapted coordinate system and also of the local logarithmic generator $\eta$ that we consider.
\end{rem}
The \emph{Newton polyhedra system} $\mathcal{N}_{M,E;\mathcal{F}}$ is the family $\{N_J\}_{J\in \mathcal{H}_{M,E}}$.
\begin{rem}
The construction of the $N_J$ is compatible with the natural projections $\text{pr}:\mathbb{R}^{J'}\to \mathbb{R}^{J}$, when $J \subset J'$, in the sense that $N_{J}=\text{pr}(N_{J'})$. This property is the essential condition asked in the general theory of polyhedra systems in \cite{Mol}.
\end{rem}

The non-degeneracy conditions for the definition of Newton non-degenerate foliated spaces, concern to the ``weighted initial forms'' associated to the compact faces of the polyhedra. Next, we give the precise statements and definitions.

\subsection{Weighted Initial Forms}
Given $J\subset I$, a \emph{$J$-vector of weights} is a linear map $\rho:\mathbb{R}^J \to \mathbb{R}$ such that $\rho(\sigma) > 0$, for every  non-zero $\sigma \in \mathbb{R}_{\geq 0}^J$. Note that the \emph{set of values} $\mathcal{V}_\rho=\rho(\mathbb{Z}^J_{\geq 0})$ has the same ordinal as $\mathbb{Z}_{\geq 0}$. 

Let us fix  $J\in \mathcal{H}_{M,E}$ and a $J$-vector of weights $\rho:\mathbb{R}^J \to \mathbb{R}$. We build the ``espaces étalé''  $\mathfrak{G}^{\rho}\to S_J$ and  $\mathfrak{A}^\rho\to S_J$ of $\rho$-weithed initial forms, as follows.

We define the fibers at a point $p\in S_J$. Given $v\in \mathcal{V}_\rho$, we consider $I_v \subset \mathcal{O}_{M,p}$ the ideal generated by the monomials $\mathbf{x}^{\sigma}$ with $\rho(\sigma)\geq v$, where $(\mathbf{x},\mathbf{y})$ is a coordinate system at $p$ adapted to $E$. On this way, we obtain a $\rho$-weighted filtration $\{I_v\}_{v\in \mathcal{V}_\rho}$ of the local ring $\mathcal{O}_{M,p}$. The associated $\rho$-graded algebra is
$$
\mathfrak{G}^{\rho}_{p}=\oplus_{v\in \mathcal{V}_\rho} \mathfrak{G}^{\rho}_{p}(v),
$$
where $\mathfrak{G}^{\rho}_{p}(v)=I_{v}/I^+_{v}=I_v/I_{v^+}$, with $v^+=\min\{b\in \mathcal{V}_{\rho}; b > v\}$. Note that we have $\mathfrak{G}^{\rho}_{p}(0)\simeq \mathcal{O}_{S_J,p}$, hence $\mathfrak{G}^{\rho}_{p}$ is an $\mathcal{O}_{S_J,p}$-graded algebra. Moreover, there is an isomorphism
\begin{equation}\label{eq:isomorfismos}
\mathfrak{G}^{\rho}_{p} \simeq \mathcal{O}_{S_J,p}[\mathbf{T}], \quad \mathbf{T}=(T_j)_{j\in J}
\end{equation}
in the category of $\mathcal{O}_{S_J,p}$-graded algebras, where the class $\mathbf{x}^\sigma+I_{\rho(\sigma)}^+\in \mathfrak{G}^{\rho}_p(\rho(\sigma))$ is sent to the monomial $\mathbf{T}^\sigma$ with weight $\rho(\sigma)$. Thanks to the isomorphisms in Equation \ref{eq:isomorfismos}, the disjoint union 
$$
\mathfrak{G}^{\rho}=\cup_{p\in S_J}\{p\}\times \mathfrak{G}^{\rho}_{p}
$$
has a topology such that the natural projection $\mathfrak{G}^{\rho} \to S_J$ is an ``espace étalé''. We obtain that $\mathfrak{G}^{\rho}$ is an $\mathcal{O}_{S_J}$-graded algebra locally isomorphic to $\mathcal{O}_{S_J}[\mathbf{T}]$, with the weights given by $\rho$. 

In a similar way, we build the $\mathfrak{G}_p^\rho$-graded module $\mathfrak{A}_p^{\rho}=\oplus_v \mathfrak{A}_p^{\rho}(v)$, obtained from the free $\mathcal{O}_{M,p}$-module $\Omega_{M,p}^1 (\log E)$ of the germs at $p$ of logarithmic one-forms with poles along $E$. As before, we globalize this construction to obtain a $\mathfrak{G}^{\rho}$-graded module $\mathfrak{A}^\rho$ with fibers $\mathfrak{A}^\rho_p$. For more details see \cite{Mol2}.
\begin{rem}
	The constructions of $\mathfrak{G}^{\rho}$ and $\mathfrak{A}^{\rho}$ do not depend on the choices of local coordinate systems adapted to $E$.
\end{rem}
For a non-zero logarithmic germ of one-form $\eta\in \Omega_{M,p}(\log E)$, the \emph{$\rho$-value $\nu_\rho(\eta)$ of $\eta$} is the maximum of the values $v$ in $\mathcal{V}_\rho$ such that $\eta \in I_v\Omega_{M,p}(\log E)$. Then, there is a well-defined $\rho$-initial form $L^\rho_p \eta \in \mathfrak{A}_p^{\rho}(\nu)$, where $\nu=\nu_\rho(\eta)$. By means of the isomorphisms 
$$
\mathfrak{A}_p^{\rho}\simeq (\mathfrak{G}_p^{\rho})^n\simeq (\mathcal{O}_{S_J,p}[\mathbf{T}])^n,
$$
and taking notations as in Equations \ref{eq:logaritmico} and \ref{eq:decompsigma}, we associate to $L^\rho_p\eta$ the family $\mathcal{L}^\rho_p(\eta)=(A_j[\mathbf{T}])_{j\in J\cup c(J)}$ of $\rho$-homogeneous polynomials $A_j[\mathbf{T}]$ defined by
\begin{equation}{\label{eq:formainicial}}
A_j[\mathbf{T}]=\sum\nolimits_{\sigma\in \Delta_\rho \cap \mathbb{Z}^J_{\geq 0}}a_{j,\sigma}(\mathbf{y})\mathbf{T}^\sigma,
\end{equation}
where $\Delta_\rho=\{\sigma \in \mathbb{R}^J; \; \rho(\sigma)=\nu_\rho(\eta)\}$. In this way, we can consider the set $\mathcal{L}^\rho_p(\eta)=0$ of common zeros of the $A_j[\mathbf{T}]$ as a subset of $(S_J,p)\times \mathbb{C}^J$.

Let us denote by $\mathcal{W}_{J}$ the set of $J$-vectors of weights. The compact faces of $N_J$ are precisely the sets $F_\rho=\Delta_\rho\cap N_J$, where $\rho\in \mathcal{W}_{J}$. Given $\rho, \rho'\in \mathcal{W}_{J}$, we have that $\mathcal{L}^\rho_p(\eta)=\mathcal{L}^{\rho'}_p(\eta)$ if and only if $F_\rho=F_{\rho'}$. Denote by $\mathcal{W}_{J,F}$ the set of $J$-vectors of weights $\rho$ such that $F_\rho=F$. In this way, we obtain a partition $\{\mathcal{W}_{J,F}\}$ of $\mathcal{W}_{J}$, by the compact faces $F$ of $N_J$. This gives sense to the expression ``initial form of $\eta$ with respect to a compact face $F$''.

The \emph{$\rho$-initial form $L_\rho \mathcal{F}$ of the foliation $\mathcal{F}$} is defined as the $\mathfrak{G}^{\rho}$-submodule of $\mathfrak{A}^{\rho}$ locally generated by the weighted initial forms $L^\rho_p\eta$, where $\eta$ are local logarithmic generators of $(M,E;\mathcal{F})$. 

\subsection{Non-Degeneracy Conditions}

We say that a foliated space $(M,E;\mathcal{F})$ is \emph{non-degenerate at $p\in S_J$ with respect to $\rho\in \mathcal{W}_J$} when we have that
\begin{equation}\label{eq:ND}
(\mathcal{L}^\rho_p(\eta)=0) \cap (S_J,p)\times (\mathbb{C}^*)^J=\emptyset,
\end{equation}
where $\eta$ is a local logarithmic generator of $\mathcal{F}$ as in Equation \ref{eq:logaritmico} and $\mathbb{C}^*=\mathbb{C}\setminus\{0\}$. 
Note that if $F_\rho=F_{\rho'}$, the condition also holds for $\rho'$. Thus, we say that the foliated space is \emph{non-degenerate at $p\in S_J$ with respect to a compact face $F$ of $N_J$} if it is so for the $\rho \in \mathcal{W}_{J,F}$.

\begin{rem} \label{rem:saturated}
Note that Equation \ref{eq:ND} holds for a local logarithmic generator $\eta$ of $\mathcal{F}$ if and only if it does for $\mathbf{x}^\sigma\eta$, since ${\mathcal L}^\rho_p(\mathbf{x}^\sigma\eta)=\mathbf{T}^\sigma {\mathcal L}^\rho_p(\eta)$.
\end{rem}

\begin{defin} 
A foliated space $(M,E;\mathcal{F})$ is \emph{Newton non-degenerate at a point $p\in M$} if it is non-degenerate at $p$ with respect to each compact face of $N_J$, where $J$ is such that $p\in S_J$. It is \emph{Newton non-degenerate} if it is so at every point of $M$.
\end{defin}
\begin{rem}
	Being Newton non-degenerate at a point is an open property.
\end{rem}

\section{Blowing-ups and Blowing-downs}
In this section we see that the property of being Newton non-degenerate is stable by combinatorial blowing-ups and blowing-downs. This is one of the keys for the proof of the equivalence statement.

\begin{prop} \label{teo:estabilidadNND}
Let $\pi:(M',E';\mathcal{F}')\to (M,E; \mathcal{F})$ be a combinatorial blowing-up between foliated spaces. We have that $(M,E;\mathcal{F})$ is Newton non-degenerate if and only if $(M',E';\mathcal{F}')$ is Newton non-degenerate.
\end{prop}

Write, for short, $\mathcal{H}=\mathcal{H}_{M,E}$ and $\mathcal{H}'=\mathcal{H}_{M',E'}$. Assume that the center of $\pi$ is $E_J$, with $J\in \mathcal{H}$. The blowing-ups of support fabrics is an abstract procedure introduced in \cite{Mol}, that is compatible, in a natural way, with the geometrical blowing-up of ambient spaces. Following \cite{Mol}, we have that
$$
\mathcal{H}=\mathcal{H}_s \cup \mathcal{K}_J,
$$
where $\mathcal{K}_J=\{K \in \mathcal{H}; \; J \subset K\}$ and $\mathcal{H}_s=\mathcal{H} \setminus \mathcal{K}_J$. Recall that $E_J$ is the adherence of the stratum $S_J$ and it is the union of the strata given by $E_J=\cup_{K\in \mathcal{K}_J}S_K$. The strata $S_K$, with $K\in \mathcal{K}_J$ are the ones disappearing after the blowing-up. Given $K\in \mathcal{K}_J$, the inverse image $\pi^{-1}(S_K)$ is a union of strata $S'_{K'(A)}$, with $K'(A)\in \mathcal{H}'^K_{\infty}$, where
$$
\mathcal{H}'^K_{\infty}=\{(K'(A); \; A \subsetneq J\}; \quad K'(A)=(K\setminus J) \cup A \cup \{\infty\}.
$$
Moreover $\pi(S'_{K'(A)})=S_K$, for each $K'(A)\in \mathcal{H}'^K_{\infty}$. In this way, we have that
$$
\mathcal{H}'=\mathcal{H}_s \cup \Big(\bigcup_{K \in \mathcal{K}_J} \mathcal{H}'^K_{\infty}\Big).
$$
Note that if $K\in \mathcal{H}_s$, the stratum $S_K$ has not been modified by the blowing-up and we can identify $S'_K=\pi^{-1}(S_K)$ with $S_K$; moreover, in this case, we have that $E'_{K}$ is the strict transform of $E_K$. 

Proving Proposition \ref{teo:estabilidadNND} is equivalent to show that $(M,E;\mathcal{F})$ is Newton non-degenerate at a point $p\in M$ if and only if $(M',E';\mathcal{F}')$ is so at each point $q$ of the inverse image $\pi^{-1}(p)$. We proceed in this way.

Let us take a point $p\in M$. If $p\notin E_J$, there is a unique point over $p$ and the blowing-up induces a local isomorphism; hence, we are done. 

We assume now that $p\in E_J$. There is a unique $K \in \mathcal{K}_J$ such that $p$ belongs to the stratum $S_K$. The inverse image $\pi^{-1}(p)$ intersects all the strata $S'_{K'(A)}$, with $A \subsetneq J$. Let us give a partition $\{\mathcal{W}_K^A\}_{A\subsetneq J}$ of $\mathcal{W}_K$ and bijections 
$$
\phi_K^A: \mathcal{W}_K^A\rightarrow \mathcal{W}_{K'(A)}, \quad A\subsetneq J
$$
such that the following ``stability'' property holds:
\begin{quote}
	\textbf{(S)} Let us consider a $K$-vector of weights $\rho \in \mathcal{W}_K^A$ and its image $\rho'=\phi_K^A(\rho)$. The foliated space $(M,E;\mathcal{F})$ is non-degenerate at $p$ with respect to $\rho$ if and only if $(M',E';\mathcal{F}')$ is non-degenerate at each $q \in \pi^{-1}(p)\cap S'_{K'(A)}$ with respect to $\rho'$.
\end{quote}
Once this is achieved, we end the proof of Proposition \ref{teo:estabilidadNND} as follows.

$\bullet $ Assume that $(M,E; \mathcal{F})$ is Newton non-degenerate at $p$ and take $q\in \pi^{-1}(p)$; let us see that the transformed foliated space $(M',E';\mathcal{F}')$ is also Newton non-degenerate at $q$. There is $A\subsetneq J$ such that $q$ belongs to the stratum $S'_{K'(A)}$. We need to see that $(M',E';\mathcal{F}')$ is non-degenerate at the point $q$ with respect to any $K'(A)$-vector of weights $\rho'$. Taking $\rho=(\phi_K^A)^{-1}(\rho')$, we know that $(M,E; \mathcal{F})$ is non-degenerate at $p$ with respect to $\rho$. Now by property \textbf{``S''} we are done.

$\bullet$ Conversely, suppose that $(M',E';\mathcal{F}')$ is Newton non-degenerate at each point $q\in \pi^{-1}(p)$; let us see that $(M,E; \mathcal{F})$ is Newton non-degenerate at $p$. We have to see that $(M,E; \mathcal{F})$ is non-degenerate at $p$ with respect to any $K$-vector of weights $\rho$. In view of the partition of $\mathcal{W}_K$, there is $A\subsetneq J$ such that $\rho\in \mathcal{W}_K^A$. Taking $\rho'=\phi_K^A(\rho)$, we have that $(M',E';\mathcal{F}')$ is non-degenerate at each $q\in \pi^{-1}(p)\cap S'_{K'(A)}$ with respect to $\rho'$. By property \textbf{``S''} we conclude.

Below we show the existence of a partition $\{\mathcal{W}_K^A\}_{A\subsetneq J}$ of $\mathcal{W}_K$ and bijections $\phi_K^A$ with the property ``\textbf{S}''.

For each $A \subsetneq J$, we define the subset $\mathcal{W}_K^A \subset \mathcal{W}_K$ as being the set of $K$-vectors of weights $\rho \in \mathcal{W}_K$ such that there is a number $r_\rho >0$ with the properties:
$$
\rho(\sigma_{K,j})=r_\rho, \;  j \in J\setminus A; \qquad \rho(\sigma_{K,j})>r_\rho, \; j \in A,
$$
where $\{\sigma_{K,j}\}_{j\in K}$ is the standard basis of the $\mathbb{R}$-vector space $\mathbb{R}^K$. The bijection $\phi_K^A: \mathcal{W}_K^A\rightarrow \mathcal{W}_{K'(A)}$ is given by $\phi_K^A(\rho)=\rho'$, where
$$
\rho'(\sigma_{K'(A),j})=\left\{
\begin{array}{cl}
\rho(\sigma_{K,j}), & j \in K \setminus J \\
\rho(\sigma_{K,j})-r_\rho, & j \in A \\
r_\rho, & j=\infty.
\end{array}\right.
$$
Let us prove that Property ``\textbf{S}'' holds. Recall that we have fixed a $K$-vector of weights $\rho\in \mathcal{W}_K^A$ and its image $\rho'=\phi_K^A(\rho)\in \mathcal{W}_{K'(A)}$.

Now we describe the relationship between the initial forms $L_\rho\mathcal{F}$ and $L_{\rho'}\mathcal{F}'$, working in local coordinates at $p$ and at $q$. We consider a local logarithmic generator $\eta$ of $\mathcal{F}$ at $p$; the foliation $\mathcal{F}'$ is defined by the pull-back $\eta'=\pi^*\eta$, more precisely, a local logarithmic generator of $\mathcal{F}'$ at $q$ is given by dividing $\eta'$ by a power of a local equation of the exceptional divisor (recall that $E'$ contains the exceptional divisor of $\pi$). Note that, in view of Remark \ref{rem:saturated}, we can work directly with $\eta'$, in order to deal with the non-degeneracy condition of $L_{\rho'}\mathcal{F}'$.

We write from now $K'=K'(A)$. Let $(\mathbf{x},\mathbf{y})$ be a coordinate system adapted to $E$ at $p$. Let us describe the morphism $\pi$ locally at $q\in \pi^{-1}(p)\cap S_{K'}$ by means of coordinates $(\mathbf{x}',\mathbf{y}')$ adapted to $E'$. We can place $q$ at some of ``the standard charts of the blowing-up'', that are parameterized by the elements of $J\setminus A$. Then, there is $j_0\in J\setminus A$ and scalars $\lambda_j \in \mathbb{C}^*$, for $j\in J_A=J\setminus (A\cup\{j_0\})$, such that the equations of $\pi$ are the following ones:
$$
\begin{array}{lcl}
	x_{j}=x'_j, \; y_{\ell}=y_{\ell}',& \text{ for }& j\in K\setminus J \text{ and }\ell\in c(K), \\
	x_{j_0}=x_{\infty}',\; x_{j}=x_{\infty}'x_j', &\text{ for } & j\in A,  \\
	x_{j}=x_{\infty}'(y_j'+\lambda_j), &\text{ for } & j\in J_A. 
\end{array}
$$

Write $\nu=\nu_\rho(\eta)$ and $\Delta_\rho=\{\sigma \in \mathbb{Z}_{\geq 0}^K; \; \rho(\sigma)=\nu\}$. We can decompose $\eta$ as a sum $\eta=\eta_0+\tilde{\eta}$ such that $\nu_\rho(\tilde{\eta})>\nu$ and
$$
\eta_0=\sum_{\sigma\in \Delta_{\rho}}\eta_{\sigma}\mathbf{x}^{\sigma}; \quad \eta_\sigma=\sum_{j\in K}a_{j,\sigma}(\mathbf{y})\frac{dx_j}{x_j}+\sum_{j\in c(K)}a_{j,\sigma}(\mathbf{y})dy_j.
$$
Given $\sigma\in \mathbb{Z}_{\geq 0}^K$, we denote $\lambda(\sigma)$ to the element of $\mathbb{Z}_{\geq 0}^{K'}$ such that $\lambda(\sigma)(j)=\sigma(j)$, if $j\in K'\setminus \{\infty\}$ and $\lambda(\sigma)(\infty)=\sum_{j\in J}\sigma(j)$. We have that $\pi^*\mathbf{x}^\sigma=U_{\sigma}{\mathbf{x}'}^{\lambda(\sigma)}$, where
$$
U_\sigma=\prod\nolimits_{j\in J_A}(y'_j+\lambda_j)^{\sigma(j)}. \; 
$$
Note that the equality $\rho'(\lambda(\sigma))=\rho(\sigma)$ holds. Moreover, we have that $\nu_{\rho'}(\eta')=\nu_\rho(\eta)=\nu$; more precisely, if we write $\eta'=\pi^*\eta_0+\pi^*\tilde{\eta}$, we see that
$$
\nu_{\rho'}(\pi^*\eta_0)=\nu; \quad \nu_{\rho'}(\pi^*\tilde{\eta})>\nu.
$$
Hence, the $\rho'$-initial form of $\eta'$ coincides with the $\rho'$-initial form of $\pi^*\eta_0$. The pull-back of $\eta_0$ by $\pi$ is given by 
$$
\pi^*\eta_0=\sum_{\sigma\in \Delta_{\rho}}\pi^*\eta_{\sigma}U_\sigma{\mathbf{x}}^{\lambda(\sigma)}.
$$
We denote $\Delta'_{\rho'}=\{\sigma' \in \mathbb{Z}_{\geq 0}^{K'}; \; \rho'(\sigma')=\nu\}$. Noting that $\lambda$ defines a bijection $\lambda:\Delta_\rho\to\Delta'_{\rho'}$, we can write
$$
\pi^*\eta_0=\sum_{\sigma'\in \Delta'_{\rho'}}\bar{\eta}_{\sigma'}{\mathbf{x}'}^{\sigma'},
$$
where 
$\bar{\eta}_{\sigma'}=U_\sigma\pi^*\eta_\sigma$, if $\sigma'=\lambda(\sigma)$;
the expression of $\bar{\eta}_{\sigma'}$ is given by
\begin{eqnarray} \label{eq:blow-up}
U_\sigma^{-1}\bar{\eta}_{\sigma'} &=&
\sum_{K'\setminus\{\infty\}}a_{j,\sigma}(\mathbf{y}')\frac{dx'_j}{x'_j}+
\Big(\sum_{j\in J}a_{j,\sigma}(\mathbf{y}')\Big)\frac{dx'_{\infty}}{x'_{\infty}}+ \\
\nonumber & & +\sum_{j\in J_A} \frac{a_{j,\sigma}(\mathbf{y}')}{y_j'+\lambda_j}dy'_j+\sum_{j\in c(K)}a_{j,\sigma}(\mathbf{y}')dy'_j,
\end{eqnarray}
where we identify the set $c(K')$ with the union $c(K)\cup J_A$. In view of Equation \ref{eq:blow-up}, the relationship between the initial forms $L^\rho_p(\eta)=(A_j[\mathbf{T}])_{j\in K\cup c(K)}$ and $L^{\rho'}_q(\eta')=(A'_j[\mathbf{T}'])_{j\in K'\cup c(K')}$ is given by
\begin{equation} \label{eq:formasiniciales}
\begin{array}{ll}
A'_\infty=\sum_{j\in J}F'_j, & \\
A'_j=F'_j, & j\in K'\setminus \{\infty\} \cup c(K)\\
A'_j=(y'_j+\lambda_j)^{-1}F'_j, & j\in J_A
\end{array}
\end{equation}
where $F'_j=\sum_{\sigma\in \Delta_\rho} a_{j,\sigma}(\mathbf{y}')U_\sigma\mathbf{T'}^{\lambda(\sigma)}$. 

The relations appearing in Equation \ref{eq:formasiniciales}, allow us to complete the proof of the Stability Property ``\textbf{S}''. More precisely, let us prove that $(M,E;\mathcal{F})$ is degenerate at $p$ with respect to $\rho$ if and only if there is a point $q\in \pi^{-1}(p)\cap S'_{K'}$ such that $(M',E';\mathcal{F}')$ is degenerate at $q$, with respect to $\rho'$.

$\bullet$ Assume that the foliated space $(M,E;\mathcal{F})$ is degenerate at $p$ with respect to $\rho$. There is $\boldsymbol{\mu}\in (\mathbb{C}^*)^K$ such that $A_j|_p(\boldsymbol{\mu})=0$, for every $j\in K \cup  c(K)$. Take the point $q\in \pi^{-1}(p)\cap S_{K'}$ defined by $\lambda_j=\mu_j/\mu_{j_0}\in \mathbb{C}^*$, for $j\in J_A$, and the vector $\boldsymbol{\mu}'\in (\mathbb{C}^*)^{K'}$ given by 
$$
\mu'_j=\mu_j/\mu_{j_0}, \; j\in A; \quad \mu'_j=\mu_j, \; j\in K\setminus J; \quad \mu'_\infty=\mu_{j_0}.
$$
We have that $A'_j|_q(\boldsymbol{\mu}')=0$, for every $j\in K' \cup c(K')$. Hence $(M',E';\mathcal{F}')$ is degenerate at $q$ with respect to $\rho'$.

$\bullet$ Consider a point $q\in \pi^{-1}(p)\cap S_{K'}$ defined by $(\lambda_j)_{j\in J_A}$, with $\lambda_j\in \mathbb{C}^*$. Suppose that $(M',E'; \mathcal{F}')$ is degenerate at $q$ with respect to  $\rho'$. There is a vector $\boldsymbol{\mu}'\in (\mathbb{C}^*)^{K'}$ such that $A'_j|_q(\boldsymbol{\mu}')=0$, for every $j\in K' \cup c(K')$. We take the vector $\boldsymbol{\mu}\in (\mathbb{C}^*)^{K}$ given by
$$
\mu_{j_0}=\mu'_{\infty}, \; \mu_{j}=\mu'_{\infty}\mu'_j,\; \mu_{\ell}=\mu'_{\ell}, \; \mu_k=\mu'_{\infty}\lambda_k; \quad j\in A, \; \ell \in K\setminus J, \; k\in J_A.
$$
We have that $A_j|_p(\boldsymbol{\mu})=0$, for every $j\in K \cup c(K)$. Hence $(M,E;\mathcal{F})$ is degenerate at $p$ with respect to $\rho$.

\section{Equivalence Theorem}
The objective on this section is to complete the proof of our main statement:
\begin{teo}\label{teo:equivalencia}
	A foliated space is Newton non-degenerate if and only if it is of logarithmic toric type.
\end{teo}

We consider a foliated space $(M, E; \mathcal{F})$ and we denote by $\mathcal{N}=\{N_J\}_{J\in \mathcal{H}}$ the Newton polyhedra system associated to it. The proof of Theorem \ref{teo:equivalencia} is a consequence of the following facts:
\begin{enumerate}[leftmargin=*]
	\item \label{fact:2} Theorem \ref{teo:equivalencia} holds when $\mathcal{N}$ is desingularized in the sense that each polyhedron $N_J$ has a single vertex. More precisely, if $\mathcal{N}$ is desingularized, then $(M,E;\mathcal{F})$ is Newton non-degenerate if and only if $\text{logSing}(\mathcal{F},E)=\emptyset$. This property can be verified in a direct way by looking at each point.
%
	\item \label{fact:3} The existence of reduction of singularities for the polyhedra system $\mathcal{N}$ with admissible centers $E_J$, that is, when the polyhedron $N_J$ has more than one vertex. This result is proved in \cite[Theorem 2]{Mol}.
	\item \label{fact:4} The compatibility between admissible blowing-ups of polyhedra systems and combinatorial log-admissible blowing-ups of foliated spaces. The main remark here is that $E_J\subset \text{logSing}(\mathcal{F},E)$ if and only if $N_J$ has more than one vertex. Hence, we have that $E_J$ is an admissible center for the polyhedra system $\mathcal{N}$ if and only if it is a log-admissible center for the foliated space. Moreover, the polyhedra system $\mathcal{N}'$ associated to the transform of the foliated space by a combinatorial log-admissible blowing-up is equal to the ``abstract'' transformed polyhedra system of $\mathcal{N}$ introduced in \cite{Mol}.

	\item \label{fact:1} The stability of being Newton non-degenerate under combinatorial blowing-ups and blowing-downs, stated in Proposition \ref{teo:estabilidadNND}.
\end{enumerate}

Let us conclude with the proof of Theorem \ref{teo:equivalencia}.

Assume first that $(M,E;\mathcal{F})$ is of logarithmic toric type. Fix a log-admissible combinatorial reduction of singularities 
$$
(M',E';\mathcal{F}') \rightarrow \cdots \rightarrow (M_1,E^1;\mathcal{F}_1) \rightarrow (M,E;\mathcal{F})
$$
given as in Equation \ref{eq:logadmisible}. Since $(M',E';\mathcal{F}')$ is logarithmically desingularized, we have that $E'_{J'} \not\subset \text{logSing}(\mathcal{F}',E')$, for any $J'\in\mathcal{H}_{M',E'}$. This implies that the polyhedra system $\mathcal{N}'$ associated to $(M',E';\mathcal{F}')$ is desingularized. By Fact (\ref{fact:2}), we have that $(M',E';\mathcal{F}')$ is Newton non-degenerate and by Fact (\ref{fact:1}) we conclude that $(M,E;\mathcal{F})$ is also Newton non-degenerate.

Suppose now that the foliated space $(M,E;\mathcal{F})$ is Newton non-degenerate.  Thanks to the compatibility stated in Fact (\ref{fact:4}) and the existence of reduction of singularities of the polyhedra system $\mathcal{N}$ stated in Fact (\ref{fact:3}), we have a sequence of log-admissible combinatorial blowing-ups 
$$
(M',E';\mathcal{F}') \rightarrow (M,E; \mathcal{F}),
$$
such that the polyhedra system $\mathcal{N}'$ associated to $(M',E';\mathcal{F}')$ is desingularized. Moreover, we know that $(M',E';\mathcal{F}')$ is Newton non-degenerate, by the stability property of Fact (\ref{fact:1}). Then, we can use Fact (\ref{fact:2}) to conclude that $(M',E';\mathcal{F}')$ is desingularized, that is, $\text{logSing}(\mathcal{F}',E')=\emptyset$. As a consequence, we have found a log-admissible combinatorial reduction of singularities. Hence $(M,E; \mathcal{F})$ is of logarithmic toric type.

\end{document}